# Almost Difference Sets in Nonabelian Groups

**Jerod Michel · Qi Wang**

**Abstract** We give two new constructions of almost difference sets. The first is a generic construction of $(q^2(q+1), q(q^2-1), q(q^2-q-1), q^2-1)$ almost difference sets occuring in certain groups of order $q^2(q+1)$ ($q$ is an odd prime power) having $(\mathbb{F}_{q^2}, +)$ as an additive subgroup. The generic construction yields several infinite families of almost difference sets, many of which occur in nonabelian groups. The second construction yields $(4p, 2p+1, p, p-1)$ almost difference sets occuring in dihedral groups of order $4p$ where $p \equiv 3 \pmod 4$ is a prime, which is interesting in the sense that there are conjectured to be no difference sets in dihedral groups (see Leung and Schmidt's work "Asymptotic Nonexistence of Difference Sets in Dihedral Groups" [J. Comb. Theory, Ser. A 99(2): 261-280(2002)]). Moreover, it turns out that some of the infinite families of almost difference sets obtained produce Cayley graphs which are Ramanujan graphs.

## 1 Introduction

Almost difference sets are an interesting subject of combinatorics and have applications in CDMA communications, error correcting codes, cryptogra-

Jerod Michel
Department of Computer Science and Engineering, Southern University of Science and Technology, Shenzhen 518055, China.
E-mail: michelj@sustc.edu.cn

Qi Wang
Department of Computer Science and Engineering, Southern University of Science and Technology, Shenzhen 518055, China.
E-mail: wangqi@sustc.edu.cn

phy and sequences. In cryptography, almost difference sets can be used to construct functions with optimal or near optimal nonlinearity [10]. In coding theory, certain almost difference sets yield cyclic codes with good parameters in the sense that they are optimal or have relatively large minimum distance [10]. Cyclic almost difference sets have also been used extensively to construct sequences with optimal autocorrelation which are used in CDMA communications [6], [11], [18], [26].

In the literature there are two different definitions of almost difference sets [7], [9]. A unification was given in [11].

**Definition 1** Let $G$ be a group of order $v$ and let $D$ be a $k$-subset of $G$. We say that $D$ is a $(v, k, \lambda, t)$ *almost difference set* if the multiset $\{* \, xy^{-1} \mid x, y \in D, x \neq y \, *\}$ contains $t$ members of $G$ with multiplicity $\lambda$ and the remaining $v - t - 1$ members of $G$ with multiplicity $\lambda + 1$.

We say that an almost difference set $D$ is *abelian* or *cyclic* (resp. *nonabelian* or *non-cyclic*) if the group $G$ is abelian or cyclic (resp. nonabelian or non-cyclic). Note also that almost difference sets can be viewed as generalizations of difference sets, i.e., a difference set can be viewed as an almost difference set with $t = 0$ or $t = v - 1$. The complement $G \setminus D$ of a $(v, k, \lambda, t)$ almost difference set is a $(v, v-k, v-2k+\lambda, t)$ almost difference set. A simple combinatorial constraint which can be used to show nonexistence of almost difference sets with certain parameters is that $(v-1)(\lambda+1) - t = k(k-1)$ must hold for any $(v, k, \lambda, t)$ almost difference set.

A great deal of known construction methods exist for almost difference sets [1], [2], [10], [11], [12], [28]. Few known construction methods, however, can be applied to nonabelian groups. To the best of our knowledge, the only known construction methods which yield infinite families of almost difference sets in nonabelian groups can be found in [7], and in Clayton's [5] recent paper in which it is shown how to construct almost difference sets in certain nonabelian groups of order $4n$ where $n \equiv 3 \pmod 4$.

Skew Hadamard difference sets can be useful in constructing almost difference sets. Ding, Wang and Pott, in [12], constructed almost difference sets in products of finite fields using skew Hadamard difference sets. Constructions of skew Hadamard difference sets can be found in [13], [14], [15], [16], [17], [23], [25], [27]. For a good survey on skew Hadamard difference sets, the reader is referred to [17].

Group ring notation is sometimes more convenient for discussing almost difference sets. For any finite group $G$, the group ring $\mathbb{Z}[G]$ is defined as the set of all formal sums of elements of $G$, with coefficients in $\mathbb{Z}$. The operations "$+$" and "$\cdot$" on $\mathbb{Z}[G]$ are given by

$$\sum_{g \in G} a_g g + \sum_{g \in G} b_g g = \sum_{g \in G} (a_g + b_g) g$$

and
$$\sum_{g \in G} a_g g \sum_{g \in G} b_g g = \sum_{g \in G} (\sum_{h \in G} a_h b_{h^{-1}g}) g.$$

where $a_g, b_g \in \mathbb{Z}$.

The group ring $\mathbb{Z}[G]$ is a ring with multiplicative identity $\mathbf{1}_G = 1 \cdot Id_G$, where $Id_G$ is the identity element of $G$, and for any subset $X \subset G$, we use the standard abuse of notation and also denote by $X$ the sum $\sum_{x \in X} x$, and by $X^{(-1)}$ the sum $\sum_{x \in X} x^{-1}$. Note that the $k$-subset $D$ of $G$ is a $(v, k, \lambda, t)$ almost difference set in $G$ whenever $DD^{-1} = k\mathbf{1}_G + \lambda T + (\lambda+1)(G - \mathbf{1}_G - T)$ for some $T \in \mathbb{Z}[G]$ with $|T| = t$.

## 2 A Generic Construction of Almost Difference Sets

In this section we give a generic construction of almost difference sets. The construction is, in a sense, a variation of the McFarland difference sets discussed in [21], which takes unions of cosets of the hyperplanes of $\mathbb{F}_q$ in certain groups that contain $(\mathbb{F}_q, +)$ as a subgroup.

**Theorem 1** *Let $q$ be a prime power and $E$ the additive group of $\mathbb{F}_{q^2}$. Let $H_1, ..., H_{q+1}$ be the subgroups (hyperplanes) of $E$ of order $q$. Let $K = \{k_1, \cdots, k_{q+1}\}$ be a group of order $q+1$ and $\tau : K \to (\mathbb{F}_q^*, \cdot)$ a group homomorphism. Define $G = E \rtimes K$ as follows: $(v, k) \cdot (w, k') = (v + \tau(k)w, k + k')$. Then*

$$D = \bigcup_{i=1}^{q+1} (E \setminus H_i) \cdot k_i \qquad (1)$$

*is a $(q^2(q+1), q(q^2-1), q(q^2-q-1), q^2-1)$ almost difference set in $G$.*

*Proof* We will use the group ring notation introduced in Section 1. We first establish some equalities. It is clear that

$$H_1 + \cdots + H_{q+1} = q\mathbf{1}_E + E. \qquad (2)$$

It is easy to see that that $H_i^2 = qH_i$, $H_iE = qE$, and $H_iH_j = E$, for $i \neq j$, hold by simple group theoretic arguments. Then we have that

$$\begin{aligned}(E - H_i)(E - H_j) &= E^2 - H_iE - EH_j + H_iH_j \\ &= q^2E - 2qE + E \\ &= (q-1)^2 E \end{aligned} \qquad (3)$$

whenever $i \neq j$ and, similarly, we can show that

$$(E - H_i)(E - H_i) = q(q-2)E + qH_i. \qquad (4)$$

Again using group ring notation, the set $D$ given by (1) becomes $D = \sum_{i=1}^{q+1}(E - H_i)k_i$. Note that, since $\tau(K)$ is contained in the kernel $\mathbb{F}_q$ of the regular spread $\{H_1, ..., H_{q+1}\}$, it follows that $k_i H_j = H_j k_i$ for all $i, j, 1 \leq i, j \leq q+1$. Thus we have

$$DD^{(-1)} = \sum_{i=1}^{q+1}(E - H_i)(E - H_i)\mathbf{1}_K + \sum_{\substack{1 \leq i,j \leq q+1 \\ i \neq j}} (E - H_i)(E - H_j)k_i k_j^{-1}$$

$$= q\sum_{i=1}^{q+1}((q-2)E + H_i)\mathbf{1}_K + (q-1)^2 \sum_{\substack{1 \leq i,j \leq q+1 \\ i \neq j}} E k_i k_j^{-1}$$

(by (4) and (3))

$$= q(q-2)(q+1)E\mathbf{1}_K + q\sum_{i=1}^{q+1} H_i \mathbf{1}_K + (q-1)^2(q+1)E(K - \mathbf{1}_K)$$

$$= q(q-2)(q+1)E\mathbf{1}_K + q(q\mathbf{1}_E + E)\mathbf{1}_k + (q-1)^2(q+1)E(K - \mathbf{1}_K)$$

(by (2))

$$= q^2(q+1)\mathbf{1}_G + (q^3 - q^2 - q)(E - \mathbf{1}_E)\mathbf{1}_K + (q^3 - q^2 - q + 1)E(K - \mathbf{1}_K).$$

□

Note moreover that, if $D$ is defined as in Theorem 1, then its complement $\overline{D}$ is an almost difference set with parameters $(q^2(q+1), q(q+1), q, q^2 - 1)$.

We give some specific examples of infinite families.

*Example 1* (i) Let $E = \mathbb{F}_{q^2}$, $K = \mathbb{Z}/(q+1)\mathbb{Z}$, and $\tau : K \to (\mathbb{F}_q^*, \cdot)$ be the group homomorphism given by $k \mapsto (-1)^k$ for all $k \in K$. Define a binary operation on the set $E \times K$ by

$$(x, k)(y, k') = (x + (-1)^k y, k + k').$$

Under this binary operation, the set $E \times K$ forms the group $G = E \rtimes K$ which is the semidirect product of $E$ and $K$ with respect to the homomorphism $\tau$. Thus, by Theorem 1, the groups $G$ so defined yield an infinite family of $(q^2(q+1), q(q^2-1), q(q^2-q-1), q^2-1)$ almost difference sets. Note that, since $\tau$ is nontrivial, the groups $G$ (and hence the almost difference sets) are nonabelian.

(ii) Let $E = \mathbb{F}_{q^2}$, $K = \mathbb{Z}/(q+1)\mathbb{Z}$, and $\tau : K \to (\mathbb{F}_q^*, \cdot)$ be the group homomorphism given by $k \mapsto 1$ for all $k \in K$. Then the binary operation given by

$$(x, k)(y, k') = (x + \tau(k)y, k + k')$$

is simply componentwise addition, under which the set $G = E \times K$ forms a group that is the direct product of $E$ and $K$, and by Theorem 1, the groups $G$ so defined yield an infinite family of abelian $(q^2(q+1), q(q^2-1), q(q^2-q-1), q^2-1)$ almost difference sets.

## 3 Almost Difference Sets in Dihedral Groups

Here we give a construction of almost difference sets in dihedral groups. A difference set $D$ in a group $G$ is called skew Hadamard if $G$ is the disjoint union of $D$, $D^{(-1)}$ and $\{0\}$. Skew Hadamard difference sets must (up to complementation) have parameters $(4\lambda-1, 2\lambda-1, \lambda-1)$ for some nonegative integer $\lambda$. Note that, for $\mathbb{Z}_n$ to contain a skew Hadamard difference set, $n$ must be a prime congruent to 3 (mod 4) (see [19]).

Let $p \equiv 3 \pmod{4}$ be a prime, and let $C_0 = \overline{C}_1 \setminus \{0\} = -C_1 \subseteq \mathbb{Z}_p$ be a skew Hadamard $(4\lambda + 3, 2\lambda + 1, \lambda)$ difference set. We again use the group ring notation that was introduced in Section 1. The following equalities are easily proved using the skew Hadamard property:

$$\overline{C}_0 \overline{C}_0^{(-1)} = \overline{C}_1 \overline{C}_1^{(-1)} = 2(\lambda+1)\mathbf{1}_{\mathbb{Z}_p} + (\lambda+1)(\mathbb{Z}_p \setminus \{0\}), \tag{5}$$

$$C_0 C_0^{(-1)} = C_1 C_1^{(-1)} = 2(\lambda+1)\mathbf{1}_{\mathbb{Z}_p} + (\lambda+1)(\mathbb{Z}_p \setminus \{0\}), \tag{6}$$

$$C_0 \overline{C}_0^{(-1)} = C_1 \overline{C}_1^{(-1)} = (\lambda+1)(\mathbb{Z}_p \setminus \{0\}), \tag{7}$$

$$\overline{C}_i \overline{C}_j^{(-1)} = \mathbf{1}_{\mathbb{Z}_p} + (\lambda+1)C_i + (\lambda+2)C_j, \tag{8}$$

$$C_i \overline{C}_j^{(-1)} = (2\lambda+1)\mathbf{1}_{\mathbb{Z}_p} + \lambda C_j + (\lambda+1)C_i, \tag{9}$$

where $i, j \in \{0, 1\}, i \neq j$.

Let $D_{2p} = \langle r, s \mid r^{2p} = s^2 = (sr)^2 = 1 \rangle$ be the dihedral group of order $4p$. Let $\phi : \mathbb{Z}_{2p} \to \mathbb{Z}_2 \times \mathbb{Z}_p$ be the map given by the Chinese Remainder Theorem, and denote $\{r^a \mid a \in \phi^{-1}(\{i\} \times C_j)\}$ by $C_{ij}$ for $i, j \in \{0, 1\}$, where $r$ is the generator of the cyclic subgroup of $D_{2p}$ of order $2p$, and denote $\{r^a \mid a \in \phi^{-1}(\{i\} \times \overline{C}_j)\}$ by $C_{i\bar{j}}$. Note that $(5) - (9)$ are easily translated to equations about the $C_{ij}$'s (and $C_{i\bar{j}}$'s) via the isomorphism $\phi$, e.g.

$$\begin{aligned}C_{01} C_{1\bar{1}}^{(-1)} &= \{r^a \mid a \in \phi^{-1}(\{0\} \times C_1)\} \{r^a \mid a \in \phi^{-1}(\{1\} \times \overline{C}_1)\}^{(-1)} \\ &= \{r^a \mid a \in \phi^{-1}(\{1\} \times (C_1 \overline{C}_1^{(-1)}))\}.\end{aligned}$$

We have the following.

**Theorem 2** *Let $p \equiv 3 \pmod{4}$ be a prime, and let $C_0 = -C_1 \subseteq \mathbb{Z}_p$ be a skew Hadamard $(4\lambda + 3, 2\lambda + 1, \lambda)$ difference set. Then the set*

$$C = (C_{0\bar{1}} \cup C_{11}) \cup s(C_{0\bar{1}} \cup C_{1\bar{1}}) \subseteq D_{2p}$$

*is a $(16\lambda + 12, 8\lambda + 7, 4\lambda + 3, 4\lambda + 2)$ almost difference set.*

*Proof* We consider the following group ring equation.

$$\begin{aligned}
CC^{(-1)} &= [(C_{0\bar{1}} \cup C_{11}) \cup s(C_{0\bar{1}} \cup C_{1\bar{1}})] [(C_{0\bar{1}} \cup C_{11}) \cup s(C_{0\bar{1}} \cup C_{1\bar{1}})]^{(-1)} \\
&= (C_{0\bar{1}} \cup C_{11})(C_{0\bar{1}} \cup C_{11})^{(-1)} + (C_{0\bar{1}} \cup C_{11})(C_{0\bar{1}} \cup C_{1\bar{1}})^{(-1)} s \\
&\quad + s(C_{0\bar{1}} \cup C_{1\bar{1}})(C_{0\bar{1}} \cup C_{11})^{(-1)} + s(C_{0\bar{1}} \cup C_{1\bar{1}})(C_{0\bar{1}} \cup C_{1\bar{1}})^{(-1)} s \\
&= (C_{0\bar{1}} \cup C_{11})(C_{0\bar{1}} \cup C_{11})^{(-1)} + s(C_{0\bar{1}} \cup C_{11})^{(-1)}(C_{0\bar{1}} \cup C_{1\bar{1}}) \\
&\quad + s(C_{0\bar{1}} \cup C_{1\bar{1}})(C_{0\bar{1}} \cup C_{11})^{(-1)} + (C_{0\bar{1}} \cup C_{1\bar{1}})^{(-1)}(C_{0\bar{1}} \cup C_{1\bar{1}}) \\
&= (C_{0\bar{1}} C_{0\bar{1}}^{(-1)} + C_{0\bar{1}} C_{11}^{(-1)} + C_{11} C_{0\bar{1}}^{(-1)} + C_{11} C_{11}^{(-1)}) \\
&\quad + 2s(C_{0\bar{1}}^{(-1)} C_{0\bar{1}} + C_{0\bar{1}}^{(-1)} C_{1\bar{1}} + C_{11}^{(-1)} C_{0\bar{1}} + C_{11}^{(-1)} C_{1\bar{1}}) \\
&\quad + (C_{0\bar{1}}^{(-1)} C_{0\bar{1}} + C_{0\bar{1}}^{(-1)} C_{1\bar{1}} + C_{1\bar{1}}^{(-1)} C_{0\bar{1}} + C_{1\bar{1}}^{(-1)} C_{1\bar{1}}).
\end{aligned}$$

We can use Equations (5)-(9) to show that

$$\begin{aligned}
&C_{0\bar{1}} C_{0\bar{1}}^{(-1)} + C_{0\bar{1}} C_{11}^{(-1)} + C_{11} C_{0\bar{1}}^{(-1)} + C_{11} C_{11}^{(-1)} \\
&= (4\lambda + 3)\mathbf{1}_{\langle r \rangle} + (2\lambda + 1)\phi^{-1}(\{0\} \times (\mathbb{Z}_p \setminus \{0\})) \\
&\quad + 2(\lambda + 1)\phi^{-1}(\{1\} \times (\mathbb{Z}_p \setminus \{0\})), \\
&C_{0\bar{1}}^{(-1)} C_{0\bar{1}} + C_{0\bar{1}}^{(-1)} C_{1\bar{1}} + C_{11}^{(-1)} C_{0\bar{1}} + C_{11}^{(-1)} C_{1\bar{1}} = 2(\lambda + 1)\langle r \rangle,
\end{aligned}$$

and

$$\begin{aligned}
&C_{0\bar{1}}^{(-1)} C_{0\bar{1}} + C_{0\bar{1}}^{(-1)} C_{1\bar{1}} + C_{1\bar{1}}^{(-1)} C_{0\bar{1}} + C_{1\bar{1}}^{(-1)} C_{1\bar{1}} \\
&= 4(\lambda + 1)\mathbf{1}_{\langle r \rangle} + 4(\lambda + 1)\phi^{-1}(\{(1,0)\}) \\
&\quad + 2(\lambda + 1)\phi^{-1}(\mathbb{Z}_2 \times (\mathbb{Z}_p \setminus \{0\})).
\end{aligned}$$

Plugging these into (10) we get

$$\begin{aligned}
CC^{(-1)} &= (8\lambda + 7)\mathbf{1}_{D_{2p}} + 4(\lambda + 1)\phi^{-1}(\{1\} \times (\mathbb{Z}_p)) \\
&\quad + (4\lambda + 3)\phi^{-1}(\{0\} \times (\mathbb{Z}_p \setminus \{0\})) + 4(\lambda + 1)s\langle r \rangle.
\end{aligned}$$

□

We give the following example.

*Example 2* Let $p = 7$. Then

$$\begin{aligned}
C &= (C_{0\bar{1}} \cup C_{11}) \cup s(C_{0\bar{1}} \cup C_{1\bar{1}}) \\
&= \{1, r^2, r^3, r^4, r^5, r^8, r^{13}, s, sr, sr^2, sr^4, sr^7, sr^8, sr^9, sr^{11}\} \subseteq D_{14}
\end{aligned}$$

is a (28,15,7,6) almost difference set.

We should point out that, besides their occurrence in nonisomorphic groups, the almost difference sets obtained in this section do not have the same cardinalities as those constructed in [5] (nor do their complements) and, although they are similar, they do not have the same parameters. Also, our dihedral construction requires the skew Hadamard property whereas the almost difference sets constructed in [5] do not.

## 4 A Note on the Cayley Graphs

A graph $\Gamma = (V, E)$ consists of a vertex set $V$ with $|V| = n$, and edge set $E$, and a relation that associates with each edge a pair of vertices. A graph is finite if its vertex set and edge set are finite. The adjacency matrix $A(\Gamma)$ of the graph $\Gamma$ is the matrix with rows and columns indexed by the vertices of $\Gamma$ with the $uv$-entry equal to the number of edges from vertex $u$ to vertex $v$. We say that $\Gamma$ is simple if it has no loops or multiple edges, and we say that $\Gamma$ is $k$-regular if the valencies $\Gamma(v)$ (number of neighbors) of each vertex $v \in V$ are all equal to $k$. The eigenvalues of the graph $\Gamma$ are defined to be the eigenvalues of $A(\Gamma)$ and are denoted by $\lambda_1(\Gamma) \geq \lambda_2(\Gamma) \geq \cdots \geq \lambda_n(\Gamma)$, and the set of all eigenvalues of $\Gamma$ is referred to as the spectrum of $\Gamma$. When $\Gamma$ is undirected, $A(\Gamma)$ is symmetric whence the eigenvalues of $\Gamma$ are real.

4.1 Expander Graphs and Ramanujan Graphs

Let $\Gamma = (V, E)$ be a simple graph. For a subset $\Omega \subseteq V$, its neighborhood is defined as

$$\Gamma(\Omega) = \{v \in V \setminus \Omega \mid v \text{ is adjacent to some } u \in \Omega\}.$$

**Definition 2** A $(|V|, k, c)$-*expander* is a $k$-regular graph $\Gamma = (V, E)$ such that every subset $\Omega \subseteq V$ with $|\Omega| \leq |V|/2$ is connected by edges of $\Gamma$ to at least $c|\Omega|$ vertices outside the set $\Omega$.

The difference $k - \lambda_2(\Gamma)$ between the degree of a $k$-regular graph $\Gamma$ and its second largest eigenvalue is called the spectral gap of $\Gamma$. Interestingly, the spectral gap gives much information about the expansion properties of a graph. The following *Expander Crossing Lemma* is proven in [20].

**Lemma 1** *Let $\Gamma = (V, E)$ be a $k$-regular graph with second largest eigenvalue $\lambda_2(\Gamma)$. Let $\Gamma = \Omega_1 \cup \Omega_2$ be a partition of $\Gamma$. Then*

$$e(\Omega_1, \Omega_2) \geq \frac{(k - \lambda_2(\Gamma))|\Omega_1| \cdot |\Omega_2|}{|V|}$$

*where $e(\Omega_1, \Omega_2)$ is the number of edges between $\Omega_1$ and $\Omega_2$.*

From this lemma one gathers that, in general, the larger the spectral gap, the better are the expansion properties, or, from a more practical perspective, the larger the spectral gap of a communication network, the less bottlenecks are hidden in it.

**Definition 3** A connected $k$-regular graph is called a *Ramanujan graph* if $\lambda_i \leq 2\sqrt{k-1}$ for all eigenvalues $\lambda_i$ of $\Gamma$ with $\lambda_i \neq k$.

Ramanujan graphs are rare and have excellent expansion properties. Some examples of Ramanujan graphs include the clique, the biclique and the Petersen Graph. Murty has a survey [24] on Ramanujan graphs; a challenging and open problem is to find more infinite families of Ramanujan graphs.

4.2 Generalized Difference Sets

**Definition 4** [4] Let $D$ be a $k(>1)$-subset of an Abelian Group $G$ and $S$ a nonempty subset of $G$. Suppose that for every nontrivial member $g \in G$,

$$\mu_g = \begin{cases} \mu_1, & \text{if } g \in S, \\ \mu_2, & \text{if } g \notin S. \end{cases}$$

where $\mu_g$ is defined to be the multiplicity of $g$ in the multiset $\{* \, xy^{-1} \mid x, y \in G, x \neq y \, *\}$. Then $D$ is called a $(|G|, |S|, k, \mu_1, \mu_2)$ generalized difference set relative to $S$.

It is easy to see that if $D \subseteq G$ is a $(|G|, |S|, k, \mu_1, \mu_2)$ generalized difference set with $|\mu_1 - \mu_2| = 1$ then $D$ is an almost difference set. For an integer $t$, let $D^t$ be the image of $D$ under the map $d \mapsto d^t$. If $D^t$ is a translation of $D$ then we say $t$ is a multiplier of $D$. If $-1$ is a multiplier of $D$ then we say that $D$ is reversible.

4.3 Cayley Graphs

Let $D$ be a $(|G|, |S|, k, \mu_1, \mu_2)$ generalized difference set relative to $S$ in an Abelian group $G$. We define the Cayley graph $\mathcal{G}(G; D, S)$ as the graph with vertex set $G$ and where two points $u$ and $v$ are adjacent if and only if $uv^{-1} \in D$. The following theorem was shown in [3].

**Theorem 3** *Let $D$ be a reversible $(|G|, |S|, k, \mu_1, \mu_2)$ generalized difference set relative to $S$ in an abelian group $G$. Then all nontrivial eigenvalues of the Cayley graph $\mathcal{G}(G; D, S)$ are contained in the set*

$$\{\pm\sqrt{k - \mu_1 + (\mu_1 - \mu_2)\chi(S)} \mid \chi \text{ is a non-principal character of } G\}.$$

*If $-\mu_1 + (\mu_1 - \mu_2)|S| < k^2 - k$ then $\mathcal{G}(G; D, S)$ is connected, and if $-\mu_1 + (\mu_1 - \mu_2)|S| < 3k - 4$ then $\mathcal{G}(G; D, S)$ is a Ramanujan graph.*

We are now ready to investigate when the Cayley graphs of the almost difference sets constructed in this paper are Ramanujan graphs.

Let $p$ be an odd prime, let $E$ denote the additive group of $\mathbb{F}_{p^2}$, and let $H_1, ..., H_{p+1}$ be the subgroups of $E$ of order $p$. Let $K = \{k_1, \cdots, k_{p+1}\}$ be a group of order $p + 1$, and let $G = E \times K$ be the direct product (under componentwise addition) of $E$ and $K$. Let $D$ be defined as in Theorem 1. Notice that

$$D^{(-1)} = \bigcup_{i=1}^{p+1}(E \setminus H_i)^{-1}k_i^{-1} = \bigcup_{i=1}^{p+1}(E \setminus H_i)k_i^{-1} = \bigcup_{i=1}^{p+1}(E \setminus H_i)k_i$$

if and only if $k_i^{-1} = k_i$ for all nontrivial $k_i \in K$. Thus if $K$ is an elementary abelian 2-group (i.e., if $p$ is a Mersenne prime) then both $D$ and $G \setminus D$ are reversible. Moreover, we have

$$-\lambda + (\lambda - (\lambda + 1)) |T| = 1 - p(p+1) < 3(p^2 + p) - 4 = 3k - 4,$$

whence the following.

**Theorem 4** *Let $p$ be an odd prime, let $E$ denote the additive group of $\mathbb{F}_{p^2}$, and let $H_1, ..., H_{p+1}$ be the subgroups of $E$ of order $p$. Let $K = \{k_1, \cdots, k_{p+1}\}$ be an elementary abelian 2-group of order $p + 1$, let $G = E \times K$ (equipped with componentwise addition), and let $D$ be defined as in Theorem 1. Then $\mathcal{G}(G; G \setminus D, T)$ is a Ramanujan graph.*

Note that Catalan's conjecture, which was proven in 2002 (see [22]), implies that if $q = p^n$ for some odd prime $p$ and some $n > 1$, then $q + 1$ cannot be a power of two, whence Theorem 4 holds only for an odd prime $p$ and not a prime power with exponent greater than one.

## 5 Concluding Remarks

We constructed $(q^2(q+1), q(q+1), q, q^2 - 1)$ almost difference sets in certain groups of order $q^2(q+1)$ containing the additive group of $\mathbb{F}_{q^2}$ as a subgroup. We constructed almost difference sets occurring in dihedral groups of order $4p$ as well, where $p \equiv 3 \pmod 4$ is prime, and having parameters $(4p, 2p+1, p, p-1)$.

As far as further constructions of almost difference sets having parameters similar to those constructed in Section 2, it may be worth checking whether there exist similar constructions in groups of order $q^2(q+1)$ which do not necessarily contain $(\mathbb{F}_{q^2}, +)$ as a subgroup. Concerning further constructions of almost difference sets in nonabelian groups of order $4n$ where $n = 4\lambda + 3$, there remain many unchecked cases for the ambient group: e.g., there may also be almost difference sets with similar parameters to those of Section 3 occurring in the groups $\mathbb{Z}_2 \times D_n$ (which is isomorphic to neither $D_{2n}$ nor $\mathbb{Z}_4 \rtimes \mathbb{Z}_n$), where $D_n$ is the dihedral group of order $2n$ (Dillon constructed Hadamard difference sets in similar groups [8]). Regarding further work in constructing Ramanujan Graphs, it seems that reversible almost difference sets are not that easy to come by, and that reversible almost difference sets whose parameters allow for its Cayley graph to be a Ramanujan graph are rarer still.

## Acknowledgement

The authors would like to thank Prof. Qing Xiang for his suggestions which led to some of the work in this manuscript. We would also like to thank the